\documentclass[11pt]{amsart}

\usepackage[T1]{fontenc}
\usepackage[utf8]{inputenc}
\usepackage{geometry}
\usepackage{amsmath,amsthm,amssymb,mathtools}
\numberwithin{equation}{section}

\usepackage[expansion=false]{microtype}
\usepackage{array,booktabs}
\usepackage{enumitem}
\setlist[itemize]{noitemsep,topsep=0pt,leftmargin=2em}
\setlist[enumerate]{noitemsep,topsep=0pt,leftmargin=2em}

\usepackage[colorlinks=true,citecolor=blue,linkcolor=blue,urlcolor=blue]{hyperref}

\theoremstyle{plain}
\newtheorem{theorem}{Theorem}[section]

\theoremstyle{definition}

\newtheorem{problem}[theorem]{Problem}
\newtheorem{conjecture}[theorem]{Conjecture}

\theoremstyle{remark}

\newcommand{\R}{\mathbb{R}}

\newcommand{\Sig}{\Sigma}
\renewcommand{\keywords}[1]{\par\smallskip\noindent\textbf{Keywords:} #1}

\begin{document}

\title{
\bfseries
Rigidity and Gap Phenomena in the Sphere--Ball Correspondence\\[0.8em]
\normalfont\large
Free Boundary Submanifolds in the Unit Ball and Their Closed Counterparts in the Sphere
}
\author{Niang Chen}
\address{Faculty of Arts and Sciences, Beijing Normal University, Zhuhai 519087, China}
\email{chenniang@bnu.edu.cn}
\date{}
\maketitle

\begin{abstract}
This survey reviews a collection of parallel phenomena between free boundary submanifolds in the Euclidean unit ball and closed submanifolds in the sphere, with particular emphasis on rigidity mechanisms, pinching thresholds, and canonical models. We do not regard the two theories as a unified system in one-to-one correspondence. Rather, we emphasize that in several typical settings---including low topology, strong pinching, spectral extremality, and symmetry reduction---the free boundary condition often forces stronger rigidity in the unit ball than in the closed setting. The exposition is organized around six interconnected themes. We first contrast the failure of the spherical Bernstein problem in high dimensions with the dimension-independent rigidity of free boundary minimal disks in the unit ball. We then discuss the parallel roles played by the Clifford torus and the critical catenoid in uniqueness, Morse index, and eigenvalue characterizations. Next, we review the transition from the Lawson--Simons stable currents method to the Bochner--Hardy techniques developed for free boundary problems, summarize pinching and gap theorems driven by the second fundamental form and its traceless part, and outline the linear comparison framework between Morse index and topology in the minimal, constant mean curvature, and weighted settings. Finally, we survey existence results obtained from group actions, isoparametric foliations, and recent equivariant eigenvalue optimization, thereby illustrating both the striking analogies and the essential boundary-driven differences between the closed spherical theory and the free boundary theory in the ball.

\end{abstract}

\section{Introduction}

Within the classical landscape of Riemannian geometry and geometric analysis, the study of minimal submanifolds has always occupied a central position. As a natural extension of the theory of minimal surfaces in Euclidean space, the rigidity and classification theory for closed minimal submanifolds in compact Riemannian manifolds---especially in the standard sphere $S^{n+1}$---reached a remarkable peak in the late 1960s and early 1970s. Simons \cite{Simons1968} introduced the Laplacian of the second fundamental form, now known as the Simons identity, and uncovered the celebrated curvature gap phenomenon. Shortly thereafter, Chern, do Carmo, and Kobayashi \cite{Chern1970}, together with Lawson \cite{Lawson1969}, achieved a complete geometric classification in the critical case, establishing the Clifford torus as the most fundamental nontrivial model after totally geodesic submanifolds in spherical geometry.

In recent years, with the development of geometric variational problems---and in particular under the impetus of Fraser and Schoen's pioneering work on the Steklov eigenvalue problem \cite{FraserSchoen2011}---attention has gradually shifted toward geometric objects with boundary, namely free boundary minimal submanifolds (FBMS). Such a submanifold $\Sigma^k \subset B^{n+1}$ has vanishing mean curvature, and its boundary $\partial \Sigma$ lies on the unit sphere $\partial B$ and meets it orthogonally. These objects are not only critical points of the relative area functional, but also serve, in a certain sense, as boundary-bearing counterparts of the theory of closed minimal submanifolds in the sphere.

The purpose of this survey is to organize systematically a number of parallel phenomena between free boundary submanifolds in the unit ball and closed submanifolds in the sphere. We do not interpret the two sides as a single unified theory in one-to-one correspondence. Instead, we emphasize the following point: in settings involving low topology, strong pinching, spectral extremality, and symmetry reduction, the free boundary condition often forces stronger rigidity in the unit ball. We develop this comparative viewpoint from six interrelated directions.

First, the uniqueness theory for the simplest models exhibits a striking contrast. For the sphere $S^n$, strong rigidity in low dimensions was established by Almgren and Calabi \cite{Almgren1966,Calabi1967}, whereas counterexamples due to Hsiang and Tomter show that the spherical Bernstein problem fails in higher dimensions \cite{Hsiang1983,Tomter1987}. By contrast, Fraser and Schoen proved that free boundary minimal disks in the unit ball remain totally geodesic in every dimension \cite{FraserSchoen2015}. This shows that in the low-topology disk setting the free boundary condition can dramatically enhance rigidity, although one should not interpret this as meaning that every problem with boundary is automatically ``more rigid'' than its closed analogue.

Second, in the characterization of critical models, we describe the parallel roles of the Clifford torus in the sphere and the critical catenoid in the unit ball. Lawson's unknottedness theorem and Brendle's geometric uniqueness theorem together characterize the Clifford torus \cite{Lawson1970,Lawson1970b,Brendle2013}. On the free boundary side, the Fraser--Li conjecture \cite{FraserLi2014} asserts the analogous uniqueness of the critical catenoid; while it remains open in full generality, substantial partial results have been obtained under various symmetry, spectral, or index assumptions. Moreover, the close parallels in Morse index (Index $5$ versus Index $4$) and in first-eigenvalue characterizations---see, for example, \cite{Urbano1990,Devyver2019,SmithZhou2019} for the index computations of the Clifford torus and the critical catenoid---further confirm the critical catenoid as the basic nontrivial model in free boundary theory.

Third, concerning topological vanishing theorems, we review the evolution of analytic methods. The paradigm shifts from the classical Lawson--Simons homology vanishing theorem based on stable currents \cite{LawsonSimons1973} to modern cohomology vanishing theorems for free boundary problems based on Bochner techniques. In this framework, the main technical obstacles are the lower bound for the Weitzenb\"ock curvature term and the control of the boundary integrals introduced by the free boundary condition. Work of Cavalcante, Mendes, and Vit\'orio \cite{Cavalcante2019} and of Chen and Ge \cite{ChenGe2022} shows that these difficulties can be overcome by a delicate intrinsic and extrinsic decomposition of the Bochner operator (thereby removing the flat normal bundle assumption) together with Hardy-type inequalities that convert boundary terms into interior energy estimates. In this way one recovers rigidity statements for manifolds with boundary that are closely parallel to the spherical case.

Fourth, in the study of curvature pinching and gap phenomena, both the theorem of Ambrozio and Nunes in the minimal case \cite{AmbrozioNunes2021} and the theorem of Barbosa, Cavalcante, and Pereira in the constant-mean-curvature case \cite{Barbosa2023} follow the classical pattern
\[
\begin{aligned}
\text{Bochner formula} \;&\to\; \text{sharp pinching}\\
&\to\; \text{gap phenomenon}\\
&\to\; \text{classification of canonical models}.
\end{aligned}
\]
These works successfully transplant Simons-type spectral gap theory from closed manifolds to free boundary and weighted settings.

Fifth, for the quantitative relation between Morse index and topological complexity, we review the extrinsic pinching framework established by Ambrozio, Carlotto, and Sharp \cite{Ambrozio2018DiffGeo,Ambrozio2018MathAnn}. By isometrically embedding the ambient manifold into Euclidean space and using harmonic $1$-forms to construct families of test functions, one can relate the Morse index of a minimal or free boundary minimal hypersurface to the dimension of its (relative) homology. In the weighted $f$-minimal setting one likewise obtains index lower bounds that grow linearly with the first Betti number \cite{Impera2020Rev,Impera2020Ann}.

Sixth, we briefly review the development of existence and construction methods. Starting from the symmetry reduction of Hsiang and Lawson \cite{HsiangLawson1971} and the classification of equivariant minimal hypersurfaces under $O(m)\times O(n)$ actions, we follow the generalization to Wang's isoparametric reduction \cite{Wang1994}, and then describe how Freidin, Gulian, and McGrath \cite{Freidin2017} and Siffert and Wuzyk \cite{SiffertWuzyk2022} transferred the resulting dynamical-systems framework to the free boundary problem in the unit ball. These developments show that the ODE/phase-plane machinery developed for complete minimal hypersurfaces in Euclidean space adapts almost verbatim to compact domains once the free boundary condition is reinterpreted as an endpoint condition for the profile curve.

Accordingly, the survey proceeds along the line of canonical models and their rigidity characterizations, homology and cohomology vanishing theorems, curvature pinching and gap phenomena, the interplay between Morse index and topology, and existence constructions via symmetry reduction. Our goal is to highlight both the parallel structures and the crucial differences between closed problems in the sphere and free boundary problems in the unit ball, while emphasizing how the free boundary condition substantially strengthens rigidity in several representative settings.

The structure of the paper is as follows. Section~\ref{sec:equator_disk} compares the uniqueness of equators in the sphere with that of equatorial disks in the unit ball. Section~\ref{sec:clifford_catenoid} discusses the dual roles of the Clifford torus and the critical catenoid. Section~\ref{sec:homology_vanishing} reviews corresponding homology and cohomology vanishing theorems on the closed and free boundary sides. Section~\ref{sec:pinching_gap} summarizes curvature pinching and gap theorems. Section~\ref{sec:morseindex} surveys index estimates and their CMC and weighted extensions. Section~\ref{sec:existence_equivariant} outlines existence constructions via symmetry reduction and isoparametric foliations in Euclidean space and the unit ball.

\subsection*{Notation and scope}
Unless otherwise specified, $A$ denotes the second fundamental form, $H$ the mean curvature, and $S=|A|^2$. In the hypersurface case we write
\[
\Phi := A-Hg,
\qquad
|\Phi|^2=|A|^2-nH^2.
\]
We write $\lambda_1$ for the first nonzero eigenvalue of the Laplace operator, $\sigma_1$ for the first nonzero Steklov eigenvalue, $b_1$ for the first Betti number, and $\mathrm{Index}$ for the (free boundary) Morse index. In Sections~\ref{sec:clifford_catenoid} and \ref{sec:pinching_gap} we mainly discuss surfaces or hypersurfaces. In Section~\ref{sec:homology_vanishing}, the meanings of $H^p(M;\mathbb R)$ and relative homology/cohomology are determined by context. Section~\ref{sec:morseindex} concerns three settings: closed minimal hypersurfaces, free boundary minimal hypersurfaces, and weighted $f$-minimal hypersurfaces; the relevant statements will be explained within each context.

\section{Equators in the Sphere and Equatorial Disks in the Unit Ball: A Rigidity Contrast}\label{sec:equator_disk}

In the geometry of the unit ball $B^{n+1} \subset \mathbb{R}^{n+1}$, equatorial hypersurfaces and equatorial disks through the origin play the role of basic geometric models. Let $P$ be an $n$-dimensional linear subspace of $\mathbb{R}^{n+1}$ passing through the origin, i.e.\ $0\in P$. We then define the following two fundamental objects.

\paragraph{Equator.}
This is the intersection of $P$ with the unit sphere $\partial B=S^n$:
\[
\Sigma_{\text{eq}} := \partial B \cap P \cong S^{n-1}.
\]
As a totally geodesic submanifold of $S^n$, it is the prototypical example of a closed minimal hypersurface.

\paragraph{Equatorial disk.}
This is the intersection of $P$ with the interior of the unit ball:
\[
\Sigma_{\text{disk}} := B^{n+1} \cap P \cong B^n.
\]
Its boundary is precisely the spherical equator, namely $\partial \Sigma_{\text{disk}}=\Sigma_{\text{eq}}$. The submanifold $\Sigma_{\text{disk}}$ is totally geodesic in Euclidean space and also meets $\partial B$ orthogonally along its boundary; it is therefore the most basic model in the theory of free boundary minimal surfaces.

Although the two objects are closely related geometrically, the study of their uniqueness and rigidity---that is, whether there exist non-equatorial minimal hypersurfaces homeomorphic to a sphere, or nonflat free boundary minimal surfaces homeomorphic to a disk---reveals a sharp distinction between closed and boundary-bearing variational problems.

\subsection{The spherical Bernstein problem: failure of high-dimensional rigidity}

The uniqueness problem for minimal hypersurfaces in the sphere was formally posed by S.~S.~Chern in his 1970s ICM report and is now known as the \textbf{spherical Bernstein problem} \cite{Chern1970,Hsiang1983}:

\begin{problem}[Spherical Bernstein Problem]
Let $\Sigma^{n-1}$ be a smooth embedded minimal hypersurface in the standard Euclidean sphere $S^n(1)$. If $\Sigma$ is topologically homeomorphic to $S^{n-1}$, must $\Sigma$ be a totally geodesic equator?
\end{problem}

In the low-dimensional case $n=3$, Almgren (1966) and Calabi (1967) gave an affirmative answer. They proved that every immersed minimal $2$-sphere in $S^3$ must be totally geodesic, hence must be an equator \cite{Almgren1966,Calabi1967}. This establishes strong rigidity in low-dimensional spheres.

However, this rigidity breaks down in higher dimensions. Using methods from equivariant differential geometry, Wu-Yi Hsiang constructed equivariant minimal hypersurfaces under suitable transformation groups and thereby disproved the conjecture in higher dimensions. In particular, he produced infinitely many families of embedded minimal hyperspheres distinct from equators in dimensions $S^4$, $S^5$, $S^6$, $S^7$, $S^8$, $S^{10}$, $S^{12}$, and $S^{14}$ \cite{Hsiang1983}. Tomter subsequently proved that such nonuniqueness is ubiquitous in even dimensions.

\begin{theorem}[Tomter, 1987 \cite{Tomter1987}]
Let $S^{2m}(1)$ be the standard Euclidean sphere of dimension $2m$. Then there exists an embedded minimal $(2m-1)$-sphere in $S^{2m}(1)$ that is not an equator.
\end{theorem}

In fact, Tomter proved that for every $n\ge 4$, the spherical Bernstein problem has a negative answer.

These counterexamples show that, without additional stability assumptions, the topological condition of being a sphere is not sufficient to determine the geometry of a minimal hypersurface in $S^n$ for $n\ge 4$.

\subsection{Equatorial disks in the unit ball: dimension-independent strong rigidity}

In sharp contrast with the failure of the spherical Bernstein problem in high dimensions, free boundary minimal disks in the unit ball exhibit remarkable rigidity. Regardless of the dimension, the combination of disk topology and the free boundary condition appears to force flatness.

For the three-dimensional Euclidean ball $B^3\subset\mathbb R^3$, Nitsche (1985) proved the classical uniqueness theorem.

\begin{theorem}[Nitsche, 1985 \cite{Nitsche1985}]
Let $\Sigma$ be an immersed minimal disk in the Euclidean unit ball $B^3$. If $\Sigma$ meets $\partial B^3$ orthogonally along its boundary (the free boundary condition), then $\Sigma$ must be a totally flat equatorial disk through the center of the ball.
\end{theorem}

Nitsche's proof relies on the Hopf differential technique and derives umbilicity from the vanishing of a holomorphic quartic differential. The result was later extended to other space forms of constant curvature. Ros and Souam (1997) proved that in a geodesic ball of $S^3$ or $H^3$, any stable capillary surface of genus zero and constant contact angle must be totally umbilical \cite{RosSouam1997}. In particular, a free boundary minimal disk must be totally geodesic.

Fraser and Schoen (2015) then extended this rigidity to arbitrary codimension, obtaining the most general theorem presently available for free boundary minimal disks.

\begin{theorem}[Fraser--Schoen, 2015 \cite{FraserSchoen2015}]
Let $u: D^2 \to B^n$ be a proper branched minimal immersion of the unit disk into the Euclidean $n$-ball ($n\ge 3$). If $u(D^2)$ meets $\partial B^n$ orthogonally, then $u(D^2)$ must be a totally geodesic equatorial disk.
\end{theorem}

The proof again uses complex-analytic techniques, showing that a higher-dimensional ambient space does not destroy disk rigidity. From the viewpoint of spectral geometry, the equatorial disk is also the unique extremal metric. Fraser and Schoen (2016) showed that it uniquely realizes the maximum value $2\pi$ of the normalized first Steklov eigenvalue $\sigma_1L$ among genus-zero surfaces with one boundary component, i.e.\ the equality case of the Weinstock inequality \cite{FraserSchoen2016}.

Taken together, these results reveal a sharp contrast between equators in the sphere and equatorial disks in the unit ball. In the spherical setting, rigidity depends strongly on dimension: in low dimensions, Almgren's theorem ensures uniqueness of minimal spheres in $S^3$, whereas for $n\ge 4$ the topology of the sphere no longer determines the equatorial model, and many nontrivial embedded examples exist. In the unit ball, by contrast, the rigidity of disk-type free boundary minimal surfaces is dimension-independent: Fraser and Schoen proved that every free boundary minimal disk in every dimension must be an equatorial disk. This shows that the free boundary condition dramatically restricts the geometry in low-topology situations. From the broader viewpoint of this survey, however, the correct conclusion is not that all problems with boundary are simply ``more rigid,'' but rather that the underlying mechanisms on the two sides are genuinely different.

\section{The Clifford Torus and the Critical Catenoid: Rigidity and Characterization}\label{sec:clifford_catenoid}

In the theory of minimal surfaces, the Clifford torus in $S^3$ and the critical catenoid in $B^3$ are the most basic nontrivial models on the closed and free boundary sides, respectively. Recent work on uniqueness and rigidity for these two geometric objects has revealed a deep parallelism between them.

\subsection{The Clifford minimal torus in the sphere}

In the unit sphere $S^3\subset\mathbb R^4$, the Clifford torus $T_{\mathrm{Cl}}$ is the product manifold $S^1(1/\sqrt{2})\times S^1(1/\sqrt{2})$. More explicitly, if one views $S^3$ as the unit sphere in $\mathbb R^4$, then the Clifford torus is given by \cite{Brendle2013,Lawson1970b}
\begin{equation}
T_{\mathrm{Cl}} = \left\{ (x_1, x_2, x_3, x_4) \in S^3 : x_1^2 + x_2^2 = x_3^2 + x_4^2 = \frac{1}{2} \right\}.
\end{equation}
It is an embedded minimal surface of genus one in $S^3$, with principal curvatures $1$ and $-1$.

In 1970, H.~Blaine Lawson Jr., while studying the topology of minimally embedded surfaces in $S^3$, proved that every compact embedded minimal surface of genus $g$ in $S^3$ is topologically standard, i.e.\ unknotted \cite{Lawson1970}. In the case $g=1$, this means that every embedded minimal torus is topologically equivalent to the Clifford torus. In the same year, Lawson proposed the following well-known geometric rigidity conjecture \cite{Lawson1970b}.

\begin{conjecture}[Lawson, 1970 \cite{Lawson1970b}]
Let $T\subset S^3$ be an embedded minimal torus. Then there exists an isometry $F$ of $S^3$ such that $F(T)=T_{\mathrm{Cl}}$. In other words, up to ambient isometries, the Clifford torus is the unique embedded minimal torus in $S^3$.
\end{conjecture}

This conjecture remained one of the central problems in minimal surface theory for more than forty years, until it was completely resolved by Simon Brendle in 2013.

\begin{theorem}[Brendle, 2013 \cite{Brendle2013}]
Let $F:\Sigma\to S^3$ be an embedded minimal torus in $S^3$. Then $F(\Sigma)$ must be congruent to the Clifford torus.
\end{theorem}

Brendle's proof settled the conjecture and established the central role of the Clifford torus in geometric uniqueness. Combined with Lawson's unknottedness theorem, it makes clear the distinction between uniqueness of topological type and uniqueness of geometric model.

Before Brendle's theorem, the Clifford torus had already been characterized in several variational ways. Urbano (1990) showed that it is singled out by Morse index.

\begin{theorem}[Urbano, 1990 \cite{Urbano1990}]
Let $\Sigma$ be a compact orientable minimal surface in $S^3$ that is not totally geodesic. Then its Morse index satisfies $\mathrm{Index}(\Sigma)\ge 5$. Equality holds if and only if $\Sigma$ is the Clifford torus.
\end{theorem}

From the point of view of spectral geometry, the Clifford torus also exhibits a distinguished rigidity. Takahashi's classical theorem (1966) states that if a submanifold is minimally immersed in $S^n$, then its coordinate functions are eigenfunctions of the Laplace operator with eigenvalue equal to the dimension of the submanifold \cite{Takahashi1966}. Thus, for a minimal surface in $S^3$, the coordinate functions satisfy $\Delta x=2x$. Montiel and Ros (1986) studied whether this eigenvalue $\lambda=2$ is the first nonzero eigenvalue $\lambda_1$, a question closely related to Yau's conjecture on minimal embeddings, and proved the following uniqueness theorem for tori.

\begin{theorem}[Montiel \& Ros, 1986 \cite{MontielRos1986}]
Let $M$ be a minimal torus immersed in $S^3$. If the immersion is given by first eigenfunctions of the Laplace operator, equivalently if $\lambda_1(M)=2$, then $M$ must be the Clifford torus.
\end{theorem}

This theorem shows that among all minimal tori, the Clifford torus is uniquely characterized by the property that its immersion into the sphere is induced by first eigenfunctions. Moreover, in their celebrated proof of the Willmore conjecture, Marques and Neves (2014) showed that the Clifford torus is the unique minimizer of the Willmore energy (equivalently, of the area in $S^3$) among surfaces of genus $g\ge 1$ \cite{MarquesNeves2014}. Altogether, these results provide a coherent geometric picture of the Clifford torus as the simplest nontrivial minimal surface in $S^3$.

\subsection{The critical catenoid in the unit ball and the Fraser--Li conjecture}

In the theory of free boundary minimal surfaces in the unit ball $B^3\subset\mathbb R^3$, the \textbf{critical catenoid} $\Sigma_{\mathrm{cc}}$ plays a role closely analogous to that of the Clifford torus. It is a compact embedded \emph{annulus} obtained by scaling the Euclidean catenoid so that it meets $\partial B^3$ orthogonally. Analytically, it can be described as the image of the conformal map $X:[-t_0,t_0]\times\mathbb R\to\mathbb R^3$ given by
\begin{equation}
X(t, \theta) = a_0 \cosh(t) \cos(\theta) e_1 + a_0 \cosh(t) \sin(\theta) e_2 + a_0 t e_3,
\end{equation}
where $\{e_1,e_2,e_3\}$ is the standard orthonormal basis of $\mathbb R^3$. The parameter $t_0$ is the unique positive solution of the transcendental equation $t\sinh(t)=\cosh(t)$, and the scaling constant $a_0=(t_0\cosh(t_0))^{-1}$ ensures that the boundary lies on the unit sphere and satisfies the free boundary condition.

Using compactness methods, Fraser and Li (2014) formulated the free boundary analogue of Lawson-type uniqueness:

\begin{conjecture}[Fraser--Li, 2014 \cite{FraserLi2014}]
Up to ambient congruence, the critical catenoid is the unique properly embedded free boundary minimal annulus in $B^3$.
\end{conjecture}

This conjecture remains open in full generality. However, significant partial results have been obtained by imposing additional spectral, symmetry, or index conditions.

Fraser and Schoen linked free boundary minimal surfaces to the Steklov eigenvalue problem and proved that if the coordinate functions of an embedded free boundary minimal annulus are first Steklov eigenfunctions, equivalently if $\sigma_1=1$, then the surface must be the critical catenoid \cite{FraserSchoen2016}. Building on this, different symmetry assumptions lead to a sequence of local-to-global uniqueness theorems.

\begin{theorem}[McGrath, 2018 \cite{McGrath2018}]
If an embedded free boundary minimal annulus in $B^n$ ($n\ge 3$) is symmetric with respect to $n$ mutually orthogonal hyperplanes, then it must be the critical catenoid.
\end{theorem}

Kusner and McGrath later used an analysis of first Steklov eigenspaces to prove another representative rigidity theorem: if an embedded free boundary minimal annulus in $B^3$ is invariant under the antipodal map, then it must be the critical catenoid \cite{KusnerMcGrath2024}. Seo further weakened the symmetry assumptions by proving that if an embedded free boundary minimal annulus in $B^3$ has two congruent boundary components and is symmetric under reflection in a plane, then it is again forced to be the critical catenoid \cite{Seo2021}.

As for the variational properties of the critical catenoid, Devyver \cite{Devyver2019} and Smith and Zhou \cite{SmithZhou2019} independently computed that its Morse index is exactly $4$. On this basis, Tran obtained the following index characterization.

\begin{theorem}[Tran, 2020 \cite{Tran2020}]
The critical catenoid is the unique free boundary minimal annulus in $B^3$ with Morse index $4$.
\end{theorem}

This theorem shows that, in the unit ball $B^3$, the critical catenoid is the simplest nontrivial free boundary minimal surface model after the equatorial disk, whose index is $1$.

Unlike Lawson's conjecture in $S^3$, uniqueness fails immediately on the free boundary side once the embeddedness assumption is removed. Fern\'andez, Hauswirth, and Mira (2023) constructed a family of non-rotationally symmetric immersed free boundary minimal annuli, showing that embeddedness is indispensable in such rigidity theorems \cite{Fernandez2023}.

In summary, the Clifford torus in $S^3$ and the critical catenoid in $B^3$ play highly parallel roles as canonical models. In the closed setting, Lawson's unknottedness theorem and Brendle's uniqueness theorem resolve, respectively, the topological and geometric levels of the problem; in particular, Brendle's result gives a complete affirmative answer to Lawson's conjecture. In the free boundary setting, the Fraser--Li conjecture remains open, but a coherent body of partial results---spectral characterization via $\sigma_1=1$, symmetry-based uniqueness under various reflection or antipodal invariance assumptions, and index characterization with $\mathrm{Index}=4$---collectively point toward an affirmative resolution. At the level of Morse index, the Clifford torus is uniquely characterized by $\mathrm{Index}=5$, while the critical catenoid is uniquely characterized by $\mathrm{Index}=4$. At the spectral level, the two models are tied to extremal properties of the first nonzero Laplace and Steklov eigenvalues, respectively. This parallelism indicates that the two canonical models play strikingly similar structural roles, even though the proof techniques on the two sides are not simple translations of one another: the closed case relies more strongly on global embedding theory and the classification of minimal tori, whereas the free boundary case draws on the Steklov spectrum, reflection principles, and the orthogonality condition at the boundary.

\section{Homology and Cohomology Vanishing Theorems: From Closed Submanifolds in the Sphere to Free Boundary Submanifolds}\label{sec:homology_vanishing}

In submanifold geometry, the relation between topological structure and curvature pinching conditions forms one of the core themes of rigidity theory. The classical paradigm in this direction was established by Lawson and Simons for submanifolds of the sphere. More recent work on free boundary submanifolds in the unit ball, by introducing refined boundary analytic techniques, has produced a highly parallel theoretical architecture.

\subsection{Homology vanishing for closed spherical submanifolds and Simons-type rigidity}

For a closed submanifold $M^n$ in the Euclidean sphere $S^{n+k}$, Lawson and Simons used the theory of stable currents from geometric measure theory to establish a deep relation between the norm of the second fundamental form and integer-coefficient homology groups. The essential idea is that under suitable curvature pinching conditions, no nontrivial stable integral currents can exist, which forces the relevant homology groups to vanish.

\begin{theorem}[Lawson \& Simons, 1973 \cite{LawsonSimons1973}]
Let $M^n$ be a compact submanifold immersed in the unit sphere $S^{n+k}$. If the squared norm of its second fundamental form satisfies the pointwise inequality
\begin{equation}
\|A\|^2 < \min\{p(n-p), 2\sqrt{p(n-p)}\},
\end{equation}
where $p$ and $q$ are positive integers with $p+q=n$, then for every finitely generated abelian group $G$ one has $H_p(M;G)=H_q(M;G)=0$. In particular, if $\|A\|^2<\min\{n-1,2\sqrt{n-1}\}$, then $M^n$ is a homotopy sphere.
\end{theorem}

The classical Lawson--Simons bound is an absolute constant and does not involve the mean curvature $H$. Vlachos (2007) observed that the standard product manifolds $S^k(r)\times S^{n-k}(s)$, namely the generalized Clifford tori, provide sharper geometric models. By introducing a critical function $a(n,k,|H|,c)$ depending on the mean curvature, he improved the rigidity theorem above and showed that topological rigidity persists even when the norm of the second fundamental form is large, provided its growth is controlled by $H$.

\begin{theorem}[Vlachos, 2007 \cite{Vlachos2007}]
Let $M^n$ ($n\ge 4$) be a compact orientable submanifold of the space form $Q_c^{n+m}$ with $c\ge 0$. If the squared norm $S$ of the second fundamental form satisfies
\begin{equation}
S < a(n, k, |H|, c) := nc + \frac{n^3|H|^2}{2k(n-k)} - \frac{n|n-2k||H|}{2k(n-k)} \sqrt{n^2|H|^2 + 4k(n-k)c},
\end{equation}
where $k$ is an integer satisfying $0<k<n$, then $H_p(M;\mathbb Z)=0$ for all $p\in[k,n-k]$. In particular, when $k=2$, if the fundamental group is finite, then $M$ is homeomorphic to the sphere $S^n$.
\end{theorem}

This extension shows that the topological rigidity of nonminimal submanifolds is strongly regulated by mean curvature, and that the inequality is sharp on Clifford-type models $S^k(r)\times S^{n-k}(s)$.

More recently, Onti and Vlachos (2022) introduced Bochner techniques for the Laplacian on $p$-forms and obtained results sharper than those coming from stable currents. They not only proved vanishing of homology groups, but also gave explicit estimates for Betti numbers and a complete geometric classification in the critical case.

\begin{theorem}[Onti \& Vlachos, 2022 \cite{OntiVlachos2022}]
Let $M^n$ be immersed in a Riemannian manifold whose curvature operator is bounded below by $c$, and assume that $H^2+c\ge 0$. If the pinching condition $S\le a(n,p,H,c)$ holds, then the Betti numbers satisfy $b_i(M)\le \binom{n}{i}$. Moreover, if a compact minimal immersion into the sphere $S^m$ satisfies $S\le n$, then $M$ must belong to one of the following three classes:
\begin{enumerate}
    \item a real homology sphere;
    \item a Clifford torus $T_{p,n-p}$;
    \item the standard Veronese embedding of a projective plane (for example, $\mathbb{C}P^2\to S^7$).
\end{enumerate}
\end{theorem}

\subsection{Cohomology vanishing for free boundary submanifolds in the unit ball}

For a compact submanifold $M^n$ in the unit ball $B^{n+k}$ satisfying the free boundary condition, namely $\partial M\perp \partial B$, the study of topological rigidity relies mainly on Bochner techniques. Unlike the spherical case, the main difficulty here lies in handling the boundary integral terms that appear in the Bochner formula.

Cavalcante, Mendes, and Vit\'orio (2019) were the first to prove, under the assumption of flat normal bundle, that pinching of the traceless second fundamental form $\Phi$ forces the vanishing of de Rham cohomology groups with real coefficients. Their proof uses weighted Hardy inequalities to control the boundary terms.

\begin{theorem}[Cavalcante et al., 2019 \cite{Cavalcante2019}]
Let $M^n$ ($n\ge 3$) be a compact free boundary submanifold in $B^{n+k}$ with flat normal bundle. If the traceless second fundamental form satisfies $\|\Phi\|^2<\frac{np}{n-p}$, where $1\le p\le \lfloor n/2\rfloor$, then
\[
H^p(M;\mathbb R)=H^{n-p}(M;\mathbb R)=0.
\]
\end{theorem}

The proof is driven by the Bochner method. One begins with a harmonic $p$-form $\omega$ on $M$ satisfying absolute or relative boundary conditions and applies the Weitzenb\"ock formula with boundary terms to derive an integral identity. Two technical ingredients are decisive. First, one uses Lin's extrinsic curvature estimate under the flat normal bundle assumption to control the Weitzenb\"ock curvature term in terms of $\Phi$ and the mean curvature $H$. Second, one resolves the boundary integral difficulty created by the free boundary condition by means of the Hardy-type inequality of Batista, Mirandola, and Vit\'orio, which estimates $\int_M |\omega|^2$ by a combination of the interior gradient integral $\int_M |\nabla |\omega||^2$, the boundary contribution, and mean-curvature terms. Together with a refined Kato inequality, these estimates yield a relation of the form ``a sum of nonnegative terms is less than or equal to zero,'' forcing $\omega\equiv 0$ once $|\Phi|^2$ is below the critical threshold $\frac{np}{n-p}$. The vanishing of cohomology then follows from Hodge--de Rham theory.

Chen and Ge (2022) subsequently introduced a new Hardy-type inequality and improved estimates for the Bochner operator, thereby removing the strong flat-normal-bundle assumption and optimizing the pinching constant.

\begin{theorem}[Chen \& Ge, 2022 \cite{ChenGe2022}]
Let $M^n$ ($n\ge 2$) be a compact free boundary submanifold in $B^{n+k}$. If the traceless second fundamental form satisfies
\[
\|\Phi\|^2 \le \frac{np}{n-p},
\qquad
1\le p\le \lfloor n/2\rfloor,
\]
then
\[
H^p(M;\mathbb R)=H^{n-p}(M;\mathbb R)=0.
\]
In particular, if $\|\Phi\|^2\le \frac{n}{n-1}$, then $M$ has exactly one boundary component and all middle-dimensional cohomology groups vanish.
\end{theorem}

The key idea of \cite{ChenGe2022} is to convert the apparently uncontrollable boundary term in the Bochner formula into an object that can be compared directly with the interior curvature term. The overall structure can be summarized as
\[
\begin{aligned}
\text{Bochner--Weitzenb\"ock} &+ \text{Hardy-type inequality}\\
&+ \text{refined Kato} + \text{extrinsic curvature estimate}\\
&\Longrightarrow \omega\equiv 0.
\end{aligned}
\]
Here, the new Hardy-type inequality turns $\int_{\partial M}|\omega|^2$ into an interior energy term, while the improved algebraic estimate controls the extrinsic contribution of $B[p]$ in arbitrary codimension. This removes the flat normal bundle assumption and sharpens the pinching constant. Full details can be found in \cite{ChenGe2022}.

This result is parallel to the Lawson--Simons stable currents theorem \cite{LawsonSimons1973}: once
\[
|\Phi|^2 \le \frac{np}{n-p},
\]
the middle-dimensional cohomology groups vanish. On the one hand, Chen and Ge strengthen Cavalcante--Mendes--Vit\'orio \cite{Cavalcante2019} by allowing equality in the pinching condition. On the other hand, continuing the spirit of Vlachos \cite{Vlachos2007}, they also extend the threshold to a version depending on $|H|$, thereby obtaining corresponding vanishing criteria in the nonminimal setting.

In particular, when $|H|$ is constant and $|H|\ge 1$, the threshold function reduces to
\[
R_p(|H|)=\sqrt{\frac{np}{n-p}}\;|H|,
\]
which has the same form as in the classical minimal case, although it now applies to a broader class of constant-mean-curvature submanifolds.

This section highlights a highly parallel rigidity paradigm. In the case of \emph{closed submanifolds in the sphere}, one considers a closed embedding or immersion $\Sigma^n\subset S^{n+k}$, and Lawson--Simons use the theory of stable currents to turn an extrinsic curvature pinching condition into the exclusion of nontrivial stable integral currents. This yields the vanishing of homology groups in complementary dimensions and leads to a characteristic gap phenomenon, with pointwise control of $\|A\|^2$ as the main ingredient \cite{LawsonSimons1973}. Later work by Vlachos and others further incorporates the mean curvature $H$ into the threshold function, showing that homology vanishing in the nonminimal setting can still be forced by strong extrinsic control and that the constants are often sharp on model examples such as Clifford-type product submanifolds \cite{Vlachos2007,Vlachos2025}.

By contrast, in the case of \emph{free boundary submanifolds in the unit ball}, one considers a compact free boundary $\Sigma^n\subset B^{n+k}$, where $\partial\Sigma\subset\partial B$ and $\Sigma$ meets $\partial B$ orthogonally. The guiding strategy naturally shifts from excluding stable currents to annihilating harmonic forms via the Bochner method. One applies the Weitzenb\"ock--Bochner identity with boundary terms to harmonic $p$-forms $\omega$ satisfying absolute or relative boundary conditions, replacing the spherical target statement $H_p(\Sigma;G)=0$ by the parallel free-boundary conclusion $H^p(\Sigma;\mathbb R)=0$ together with its dual-dimensional counterpart. The \emph{point of correspondence} is that the pinching quantity is still governed by the second fundamental form: on the spherical side one mainly controls $\|A\|^2$, whereas on the free boundary side it is more natural to control the traceless part $\|\Phi\|^2$, possibly coupled with $H$. The \emph{essential difference} is the boundary integral term, typically of the form $\int_{\partial\Sigma}|\omega|^2$, which has no analogue in the closed case. This is why weighted Hardy-type inequalities are indispensable: they convert the boundary contribution into an interior energy term and thereby reconstruct, within a genuinely boundary-sensitive Bochner framework, cohomology vanishing results that are almost perfectly parallel to the Lawson--Simons theory. Cavalcante--Mendes--Vit\'orio first implemented this program under additional assumptions \cite{Cavalcante2019}; Chen--Ge then removed the strong hypotheses and optimized the pinching constants through improved Hardy inequalities and extrinsic curvature estimates, bringing the free boundary rigidity theory much closer in form and conclusion to the classical picture for closed submanifolds in the sphere \cite{ChenGe2022}.

\section{Curvature Pinching Inequalities and Gap Theorems}\label{sec:pinching_gap}

In the geometry of minimal submanifolds and constant-mean-curvature (CMC) hypersurfaces, a \emph{curvature pinching} condition typically means a pointwise upper bound, as sharp as possible, on the second fundamental form $A$ or its traceless part $\Phi$. Combined with Simons-type Bochner formulas, the maximum principle, and integral identities, such bounds compress the geometric possibilities down to a small collection of canonical models. The corresponding \emph{gap phenomenon} asserts that certain natural geometric quantities---most commonly $|A|^2$ or $|\Phi|^2$---cannot occur in a nontrivial open interval, and that only the endpoint values are realized by extreme configurations. For free boundary  submanifolds, pinching and gap problems are often modeled on the spherical theory: one first derives a pinching inequality from a Bochner-type identity and then classifies the equality case, obtaining rigidity and a forbidden interval.

\subsection{Curvature pinching and gap theorems for submanifolds in the sphere}

The study of pinching and gap phenomena lies at the heart of rigidity theory. The subject was founded by James Simons in his landmark 1968 paper. Simons introduced powerful elliptic methods for the study of minimal submanifolds and derived the Laplacian formula for the second fundamental form of a minimal submanifold, now known as the \textbf{Simons identity} \cite{Simons1968}. On this basis he proved a striking pinching theorem for minimal hypersurfaces in the sphere.

Let $M^n\subset S^{n+1}(1)$ be a closed immersed hypersurface, let $A$ denote its second fundamental form, $H$ its mean curvature, and write
\[
S:=|A|^2.
\]
When $H\equiv 0$, i.e.\ in the minimal hypersurface case, Simons obtained the fundamental elliptic identity
\begin{equation}\label{eq:SimonsBochner}
\frac12\,\Delta S \;=\; |\nabla A|^2 \;+\; nS \;-\; S^2,
\end{equation}
which can be viewed as a Bochner formula for $A$. Integrating \eqref{eq:SimonsBochner} over the closed manifold $M$ and applying integration by parts yields
\begin{equation}\label{eq:SimonsIntegral}
\int_M S(S-n)\,d\mu \;=\; \int_M |\nabla A|^2\,d\mu \;\ge\;0.
\end{equation}
Therefore, once one imposes the pointwise pinching condition $S\le n$, one has $S(S-n)\le 0$, and combining this with \eqref{eq:SimonsIntegral} forces
\[
S(S-n)\equiv 0
\quad\text{and}\quad
\nabla A\equiv 0.
\]

\begin{theorem}[Simons, 1968 \cite{Simons1968}]\label{thm:SimonsGapSphere}
Let $M^n\subset S^{n+1}$ be a compact minimal hypersurface and write $S=|A|^2$. If $S\le n$ everywhere, then either $S\equiv 0$ or $S\equiv n$. In particular, there is no compact minimal hypersurface with $0<S<n$.
\end{theorem}

Simons' work established a spectral gap phenomenon for the squared norm $S$ of the second fundamental form. The next major question was to classify geometrically the submanifolds that attain the critical upper bound. Chern, do Carmo, and Kobayashi (1970), using moving frames, and Lawson (1969), using local rigidity analysis in Riemannian geometry, arrived independently at the same rigidity theorem.

\begin{theorem}[Chern--do Carmo--Kobayashi, 1970; Lawson, 1969 \cite{Chern1970,Lawson1969}]
Let $M^n$ be a compact minimal submanifold of the unit sphere $S^{n+p}$. If the squared norm $S$ of the second fundamental form satisfies
\[
S \le \frac{n}{2 - 1/p},
\]
then $S$ is constant, and $M$ must be one of the following:
\begin{enumerate}
    \item $S\equiv 0$, in which case $M$ is totally geodesic;
    \item $S\equiv \frac{n}{2 - 1/p}$. In particular, in the hypersurface case ($p=1$), the critical value $S=n$ corresponds to scalar curvature
    \[
    R=n(n-2),
    \]
    or, in normalized form,
    \[
    \bar R=\frac{n-2}{n-1}.
    \]
    In this case, equality occurs if and only if $M$ is a \textbf{Clifford torus}:
    \[
    M_{k,n-k} = S^k\left(\sqrt{\frac{k}{n}}\right)\times S^{n-k}\left(\sqrt{\frac{n-k}{n}}\right), \quad 1\le k<n.
    \]
\end{enumerate}
\end{theorem}

This theorem not only unifies the global integral estimates of Chern et al.\ with Lawson's local rigidity analysis, but also reveals that $S=n$ is a genuine transition threshold at which the geometry jumps from a sphere to a product manifold.

Through the Bochner-type formula for $S=|A|^2$, Simons established the basic framework for curvature pinching of minimal hypersurfaces in the sphere: under compactness or suitable constancy assumptions on $S$, the maximum principle yields a typical gap phenomenon, so that $S$ becomes rigid near the critical threshold. Once the equality case $S\equiv n$ occurs, the geometric model is classified precisely by Chern--do Carmo--Kobayashi and Lawson as a Clifford-type product of spheres; in higher codimension there are also exceptional models such as the Veronese surface. Thus Simons' analytic pinching inequality together with the equality classification of Chern--do Carmo--Kobayashi and Lawson laid down the standard paradigm
\[
\text{inequality} \;\to\; \text{pinching interval} \;\to\; \text{equality model} \;\to\; \text{classification/rigidity}.
\]

It is important to note that pinching and gap phenomena are not confined to the minimal setting. For hypersurfaces with constant mean curvature, the more natural quantity is the traceless second fundamental form $\Phi$ (see the notation in the Introduction), with $|\Phi|^2=|A|^2-nH^2$. Alencar and do Carmo proved a sharp CMC gap theorem: if $|\Phi|^2$ lies below a critical threshold determined by $H$ and the dimension, then the geometry must fall into either the totally umbilical case or an explicit product model.

\begin{theorem}[Alencar--do Carmo, 1994 \cite{Alencar1994}]\label{thm:AdCGap}
Let $M^n\subset S^{n+1}(1)$ be a closed orientable hypersurface with constant mean curvature. Choose the orientation so that $H\ge 0$, and let $\Phi$ denote the traceless second fundamental form. For each $H$, define
\[
P_H(x)=x^2+\frac{n(n-2)}{\sqrt{n(n-1)}}\,H\,x-n(H^2+1),
\]
and let $B_H$ be the square of the positive root of $P_H(x)=0$ (in particular, $B_0=n$). If $|\Phi|^2\le B_H$ everywhere on $M$, then either $|\Phi|^2\equiv 0$ or $|\Phi|^2\equiv B_H$. The case $|\Phi|^2\equiv 0$ is equivalent to $M$ being totally umbilical. In the equality case $|\Phi|^2\equiv B_H$, one has a local classification: when $H=0$ one recovers the minimal Clifford model; when $H\neq 0$ and $n\ge 3$, $M$ is locally an $H(r)$-torus,
\[
S^{n-1}(r)\times S^1\!\left(\sqrt{1-r^2}\right)\subset S^{n+1}(1),
\]
whose principal curvatures are
\[
k_1=\cdots=k_{n-1}=\frac{\sqrt{1-r^2}}{r},
\qquad
k_n=-\frac{r}{\sqrt{1-r^2}},
\]
for a suitable choice of $r$.
\end{theorem}

Theorems~\ref{thm:SimonsGapSphere} and \ref{thm:AdCGap} together display the most classical closed loop in spherical geometry:
\[
\begin{aligned}
\text{Bochner-type formula} \;&\to\; \text{sharp pinching}\\
&\to\; \text{gap phenomenon}\\
&\to\; \text{classification of equality cases}.
\end{aligned}
\]
Once the curvature quantity is pinched below the critical threshold, the geometry is forced to collapse to a very small list of models: totally geodesic, totally umbilical, Clifford, or their CMC product analogues. It is precisely this paradigm that provides the clearest and most portable prototype for the free boundary questions in the unit ball, such as whether sufficiently small curvature forces a surface to be either a flat disk or a critical catenoid, and whether one can identify forbidden intervals for area or curvature.

\subsection{Curvature pinching and gap theorems in the unit ball}

\subsubsection{The minimal case}

In close analogy with the classical gap theorem in the sphere, Ambrozio and Nunes (2021) proved a sharp pinching/gap theorem involving the second fundamental form for free boundary minimal surfaces in $B^3$ \cite{AmbrozioNunes2021}.

\begin{theorem}[Ambrozio--Nunes, 2021 \cite{AmbrozioNunes2021}]\label{thm:AmbrozioNunes}
Let $\Sig^2\subset B^3$ be a compact free boundary minimal surface, let $A$ be its second fundamental form, let $N$ be a local unit normal vector, and let $x$ denote the position vector. If
\[
|A|^2(x)\,\langle x, N(x)\rangle^2 \le 2
\]
for all $x\in \Sig$, then:
\begin{enumerate}
  \item if $|A|^2\langle x,N\rangle^2\equiv 0$, then $\Sig$ is a flat equatorial disk;
  \item if there exists $p\in \Sig$ such that $|A|^2(p)\langle p,N(p)\rangle^2=2$, then $\Sig$ is the critical catenoid.
\end{enumerate}
\end{theorem}

The Ambrozio--Nunes argument may be summarized as
\[
\begin{aligned}
|A|^2\langle x,N\rangle^2\le 2
&\Longrightarrow \mathrm{Hess}_{\Sigma}\frac{|x|^2}{2}\ge 0\\
&\Longrightarrow C(\Sigma)\text{ is totally convex}\\
&\Longrightarrow \Sigma \text{ must be either a disk or an annulus}.
\end{aligned}
\]
If $C(\Sigma)$ is a single point, then Nitsche's theorem yields an equatorial disk. If equality is attained at some point, then $C(\Sigma)$ is an interior closed geodesic, and one can use a Jacobi function, the structure of the nodal set, and unique continuation to infer rotational symmetry, which forces the surface to be the critical catenoid. The real bridge is therefore not a long local computation, but the mechanism
\[
\begin{aligned}
\text{pinching}
&\to \text{convexity of the radial function}\\
&\to \text{topological compression}\\
&\to \text{recognition of the canonical model},
\end{aligned}
\]
as emphasized in \cite{AmbrozioNunes2021,Hardt1989}.

This framework has been widely used in later work. Ambrozio and Nunes asked whether their theorem could be extended to higher dimensions or higher codimensions. The first positive answer was obtained for two-dimensional free boundary minimal surfaces in arbitrary codimension.

Barbosa and Viana extended the two-dimensional free boundary minimal surface gap theorem from $B^3$ to $B^n$ \cite{BarbosaViana2020}.

\begin{theorem}[Barbosa--Viana, 2020 \cite{BarbosaViana2020}]\label{thm:BV2020}
Let $\Sig^2\subset B^n$ be a compact free boundary minimal surface, let $A$ be its second fundamental form, and let $x^\perp$ denote the normal component of the position vector $x$. If
\[
|x^\perp|^2\,|A(x)|^2 \le 2
\]
for all $x\in\Sig$, then:
\begin{enumerate}
  \item if $|x^\perp|^2|A|^2\equiv 0$, then $\Sig$ is a flat equatorial disk;
  \item if there exists $x_0\in\Sig$ such that $|x_0^\perp|^2|A(x_0)|^2=2$, then $\Sig$ is actually contained in a three-dimensional linear subspace, and within that subspace it is congruent to the critical catenoid.
\end{enumerate}
\end{theorem}

At the level of proof architecture, this is almost a direct generalization of the three-dimensional codimension-one argument of Ambrozio--Nunes. In both cases one begins with the distance function
\[
f(x)=\frac{|x|^2}{2},
\]
converts the curvature pinching hypothesis into the convexity conclusion $\mathrm{Hess}_\Sigma f\ge 0$, and deduces that the minimum set
\[
C=\{f=\min f\}
\]
is totally convex. The free boundary condition then enforces strict boundary convexity, showing that $\Sigma$ must be either a disk or an annulus. In the disk case one invokes the uniqueness of free boundary minimal disks to conclude that $\Sigma$ is an equatorial disk.

The crucial difference is that Ambrozio--Nunes relies heavily on the low-dimensional codimension-one structure of $B^3$. They write down the eigenvalues of $\mathrm{Hess}_\Sigma f$ explicitly in terms of the principal curvatures. In the equality case, they use the cross product to construct the rotational Jacobi function
\[
u=\langle x\wedge N,e_3\rangle.
\]
Cheng's nodal set theorem then forces $u\equiv 0$, hence rotational symmetry, and therefore the critical catenoid. By contrast, Barbosa--Viana must treat higher dimensions and arbitrary codimensions in $B^n$. They therefore prove $\mathrm{Hess}_\Sigma f\ge 0$ by an algebraic inequality, no longer relying on the three-dimensional identity ``principal curvatures $=\pm |A|/\sqrt2$.'' The uniqueness of disks is replaced by the Fraser--Schoen theorem in arbitrary codimension. In the annular/equality case, the argument is changed to ``dimension reduction plus nodal-set analysis for an elliptic system'': along $C$ they construct a normal distribution $E$ and prove that it is parallel, which shows that $\Sigma$ actually lies in a fixed three-dimensional linear subspace. One then compares $\Sigma$ locally with the critical catenoid in that subspace and invokes the Hardt--Simon description of nodal sets together with unique continuation to force coincidence.

Barbosa and Viana then extended the gap phenomenon to higher-dimensional minimal submanifolds. In 2022 they obtained a topological characterization in arbitrary dimension \cite{BarbosaViana2022}.

\begin{theorem}[Barbosa--Viana, 2022 \cite{BarbosaViana2022}]\label{thm:BV2022}
Let $\Sig^k$ be a free boundary minimal submanifold of a convex domain $\Omega\subset \R^{n+1}$. Write $x^\perp$ for the normal component of the position vector and define
\[
f(x)=\frac{|x|^2}{2},
\qquad 
C(\Sig):=\Bigl\{x\in\Sig:\ f(x)=\min_{\Sig} f\Bigr\}.
\]
If
\[
|A|^2(x)\,|x^\perp|^2 \le \frac{k}{k-1}
\]
for all $x\in\Sig$, then one of the following holds:
\begin{enumerate}
  \item $C(\Sig)=\{x_0\}$ and $\Sig$ is diffeomorphic to the disk $\mathbb{D}^k$;
  \item $C(\Sig)$ is a closed geodesic and $\Sig$ is diffeomorphic to $\mathbb{S}^1\times \mathbb{D}^{k-1}$. Moreover, equality holds along $C(\Sig)$, and for $x\in C(\Sig)$ one has
  \[
  \left.\mathrm{Hess}_{\Sig} f \right|_{T_x C(\Sig)^\perp} = \frac{k}{k-1}\,\mathrm{Id}.
  \]
\end{enumerate}
In particular, if
\[
|A|^2|x^\perp|^2<\frac{k}{k-1}
\]
holds everywhere on $\Sig$, then $\Sig$ is diffeomorphic to $\mathbb{D}^k$.
\end{theorem}

Theorems~\ref{thm:BV2020} and \ref{thm:BV2022} follow the same broad strategy, both beginning with the Hessian analysis of the radial function $f(x)=|x|^2/2$ on a free boundary minimal submanifold $\Sigma^k\subset B^{n+1}$. Using the free boundary condition $\Sigma\perp \partial B$, one obtains the uniform formula
\begin{equation}
\operatorname{Hess}_{\Sigma} f(X,Y) = \langle X, Y \rangle + \langle A(X, Y), x \rangle,
\end{equation}
where $A$ is the second fundamental form. Barbosa and Viana observed that the curvature pinching conditions
\[
|x^\perp|^2|A|^2 \le 2
\qquad\text{or}\qquad
|x^\perp|^2|A|^2 \le \frac{k}{k-1}
\]
are strong enough to guarantee that $\operatorname{Hess}_{\Sigma} f$ is positive semidefinite.

This convexity implies that the minimum set
\[
C(\Sigma)=\{x\in\Sigma: f(x)=\min f\}
\]
is totally convex. As a consequence, the topology of the submanifold is compressed into two basic possibilities only: the disk $D^k$ or the solid tube $S^1\times D^{k-1}$.

Even though the overall framework is the same, the precision of the conclusions and the tools involved differ substantially. Theorem~\ref{thm:BV2020} focuses on the two-dimensional case $k=2$ but allows arbitrary codimension. After establishing the topological dichotomy ``disk versus annulus,'' the authors pursue a much stronger geometric congruence statement. In the disk case, the strict boundary convexity induced by the free boundary condition allows one to invoke Fraser--Schoen's uniqueness theorem \cite{FraserSchoen2015} and identify the surface with a flat equatorial disk. In the annulus case, the proof involves a delicate geometric construction: the authors first build a parallel normal distribution along $C(\Sigma)$ and show that $\Sigma$ is actually contained in a three-dimensional linear subspace. They then compare $\Sigma$ with the standard critical catenoid by second-order tangency and use Hardt--Simon's nodal set analysis together with unique continuation for elliptic systems \cite{Hardt1989} to conclude that the two surfaces coincide.

Theorem~\ref{thm:BV2022}, by contrast, is aimed at arbitrary dimension $k$ and therefore seeks topological rather than geometric rigidity, since in high dimensions no classification theorem comparable to Fraser--Schoen is currently available. A purely algebraic lemma, concerning inequalities between sums of squares and squares of sums, ensures that the relevant trace inequality controls the signs of all $k$ eigenvalues of $\operatorname{Hess}_{\Sigma} f$ simultaneously. Combined with Milnor's Morse theory \cite{Milnor1963}, this yields directly that $\Sigma$ is diffeomorphic either to $D^k$ or to $S^1\times D^{k-1}$. In the critical equality case, the conclusion remains essentially second-order along $C(\Sigma)$, namely
\[
\mathrm{Hess}_{\Sigma} f|_{C(\Sigma)^\perp} = \frac{k}{k-1}\mathrm{Id},
\]
rather than identifying a unique geometric model by analytic continuation as in the 2020 paper.

These conclusions may be viewed as a higher-dimensional topological version of the Ambrozio--Nunes theorem: although far more complicated minimal submanifolds may exist in high dimensions, sufficiently strong curvature pinching leaves room only for a topologically simple disk-type model or a single $S^1$-type fundamental loop. It should be emphasized, however, that the high-dimensional results currently guarantee only topological classification; the precise geometric models remain incompletely understood.

For example, in the four-dimensional ball $B^4$, one expects free boundary minimal hypersurfaces analogous to the critical catenoid, possibly with $O(2)\times O(2)$ symmetry \cite{BarbosaViana2022}. Barbosa and Viana formulated the following conjecture (for a statement, see \cite{BarbosaViana2022}):

\begin{conjecture}[Barbosa--Viana, 2022 \cite{BarbosaViana2022}]\label{conj:BV_B4}
Let $\Sig^3\subset B^4\subset \R^4$ be a capillary minimal hypersurface in the unit ball; the free boundary case corresponds to contact angle $\pi/2$. Assume that
\[
|x^\perp|^2\,|A(x)|^2 \le \frac{3}{2}
\]
for all $x\in\Sig$. Then one of the following should hold:
\begin{enumerate}
  \item $\Sig$ is congruent to the equatorial three-disk $\mathbb{D}^3$;
  \item $\Sig$ is congruent to an $O(2)\times O(2)$-invariant free boundary minimal hypersurface; such examples are topologically homeomorphic to the solid torus $\mathbb{D}^2\times \mathbb{S}^1$.
\end{enumerate}
\end{conjecture}

At present, the concrete uniqueness and existence problems in this higher-dimensional setting remain unresolved. In light of the work of Freidin, Gulian, and McGrath \cite{Freidin2017}, the main difficulty in resolving Conjecture~\ref{conj:BV_B4} is that one must distinguish a specific model from the infinitely many examples already known to exist. Freidin--Gulian--McGrath proved that in $B^4$ (that is, in the case $m=n=2$), the associated differential equation has focal singularities near the singular point, leading to oscillatory behavior of the solution curves. As a consequence, there exist infinitely many pairwise noncongruent $O(2)\times O(2)$-invariant free boundary minimal hypersurfaces $\{\Sigma_{2,2,k}\}$, all homeomorphic to solid tori \cite{Freidin2017}. Barbosa and Viana's conjecture attempts to single out a distinguished solid-torus model by the curvature condition
\[
|x^{\perp}|^2 |A|^2 \le \frac{3}{2},
\]
but in the presence of this infinite family one must prove that all higher-oscillation members violate the pinching inequality somewhere globally.

Beyond the Euclidean unit ball, analogous gap phenomena have been studied in other space forms as well. Li and Xiong (2018) extended the Ambrozio--Nunes framework to geodesic balls in hyperbolic space $\mathbb H^3$ and the hemisphere $\mathbb S^3_+$ \cite{LiXiong2018}. By introducing potential-function correction terms adapted to the ambient metric, they proved that under suitable weighted curvature pinching conditions, a free boundary minimal surface must be either a totally geodesic disk or a rotationally symmetric annulus, the latter being the natural nonflat analogue of the critical catenoid in a non-Euclidean background.

Min and Seo (2022) pushed the theory further to strictly convex domains in general Riemannian three-manifolds with upper bounds on sectional curvature, and they also included the CMC case \cite{MinSeo2022}. By analyzing pinching inequalities involving the traceless second fundamental form and the Hessian of the distance function, they obtained a unified topological rigidity theorem: under the pinching hypothesis, the surface must be homeomorphic to either a disk or an annulus. In the special case of geodesic balls in space forms, this topological result upgrades to a geometric classification, showing that the surface must be either a spherical cap or a Delaunay surface. These developments confirm that curvature-gap rigidity is not a peculiarity of Euclidean space, but rather a robust feature across several geometric backgrounds.

\subsubsection{The CMC case}

For hypersurfaces with nonzero constant mean curvature $H\neq 0$, the quantity $|A|^2$ is often no longer the most natural pinching object. Inspired by the spherical work of Alencar and do Carmo (1994) \cite{Alencar1994}, recent research has instead focused on the traceless second fundamental form $\Phi$ (see the notation in the Introduction), using $|\Phi|^2$ as the curvature pinching quantity.

For compact free boundary CMC submanifolds in the unit ball $B^{n+1}$, one frequently encounters a gap/rigidity classification of the following kind: if the traceless part $\Phi$ satisfies a suitable pointwise pinching condition, often coupled with the support function $\langle x,N\rangle$, then the hypersurface $\Sigma$ must belong to a very small explicit family of models.

Barbosa, Cavalcante, and Pereira (2023) established the basic gap theorem for CMC surfaces in the unit ball, with an inequality whose right-hand side depends explicitly on the mean curvature $H$.

\begin{theorem}[Barbosa, Cavalcante \& Pereira, 2023 \cite{Barbosa2023}]
Let $\Sigma$ be a compact free boundary CMC surface in the unit ball $B^3$, with $H\neq 0$. If the traceless second fundamental form $\Phi$ satisfies
\begin{equation} \label{eq:bcp}
\|\Phi\|^2 \langle x, N \rangle^2 \le \frac{1}{2} (2 + H \langle x, N \rangle)^2,
\end{equation}
then one of the following holds:
\begin{enumerate}
    \item either $\|\Phi\|\equiv 0$, in which case $\Sigma$ is totally umbilical, hence a spherical cap;
    \item or equality holds at some point, in which case $\Sigma$ must be part of a rotationally symmetric Delaunay surface.
\end{enumerate}
\end{theorem}

This theorem shows that a CMC surface satisfying the pinching condition can only be a spherical cap or a Delaunay surface. It may be viewed as the free boundary analogue of the Alencar--do Carmo gap theorem for closed CMC hypersurfaces in the sphere \cite{Alencar1994}.

After this foundational Euclidean result, rigidity theory quickly expanded to broader geometric settings. Andrade, Barbosa, and Pereira (2021) first extended the discussion to radially symmetric conformally Euclidean three-balls $B^3_r$ \cite{Andrade2021}. They introduced a potential function $\sigma$ adapted to the conformal metric and proved that if the traceless second fundamental form $\Phi$ satisfies a suitable weighted pinching inequality involving $\sigma$ and the support function, then a free boundary CMC surface again exhibits a dichotomy between disks and rotationally symmetric annuli. Examples in Gaussian space clarify the geometric meaning of the condition.

For ambient spaces with intrinsic curvature, Min and Seo (2022) established an analogous topological rigidity theorem for strictly convex domains in three-dimensional Riemannian manifolds. They showed that under a corresponding pinching condition on $\|\Phi\|$, a free boundary CMC surface must again degenerate topologically to either a disk or an annulus. In geodesic balls of space forms, this topological conclusion can be sharpened to a geometric classification: the surface must be either a spherical cap or a Delaunay surface \cite{MinSeo2022}.

Most recently, Freitas, Santos, and Sindeaux (2025) extended these results systematically to general rotational domains such as rotational ellipsoids \cite{Freitas2025}. They not only classified surfaces satisfying the gap condition, but also investigated corresponding higher-dimensional minimal hypersurface phenomena. More importantly, by invoking the existence of non-rotationally symmetric free boundary annuli in the Euclidean ball, they highlighted from the opposite direction the necessity of such curvature pinching assumptions for isolating the canonical geometric models.

\section{Linear Comparison Between Index and Topology: From Closed Minimal Hypersurfaces to the Weighted Free Boundary Setting}\label{sec:morseindex}

The Morse index of a minimal hypersurface, denoted $\mathrm{Index}$, measures its variational instability. In the two-sided case, $\mathrm{Index}(M)$ is the maximal dimension of a subspace of admissible normal variations on which the second variation of area is negative definite. In closed ambient manifolds with positive Ricci curvature, a conjectural picture advocated by Schoen, Marques, Neves, and others predicts that the index of a minimal hypersurface should admit an affine lower bound in terms of its topological complexity, especially the first Betti number $b_1$; see \cite{Ambrozio2018DiffGeo,Impera2020Ann} for relevant formulations and background.

Ambrozio, Carlotto, and Sharp developed a unified analytic framework based on isometric embeddings of the ambient space, significantly advancing the quantitative comparison between index and topology.

Since the original papers formulate their assumptions in terms of integral pinching inequalities involving the extrinsic geometry of the ambient manifold, we state the results below with their full hypotheses.

Unlike Ros, who in the surface case uses test functions of the form $u=\omega(V)$ with $\omega$ a harmonic $1$-form and $V$ an ambient parallel vector field, and unlike Savo, who in higher dimensions introduces more refined combinations of such functions, Ambrozio--Carlotto--Sharp systematically construct test functions from coordinates in $\Lambda^2\mathbb R^d$ under a general isometric embedding $\mathcal N^{n+1}\hookrightarrow\mathbb R^d$. If $N$ is the unit normal of $M$ and $\omega^\sharp$ is the dual vector field of a harmonic $1$-form $\omega$, they define
\[
u_{ij} \;=\; \big\langle N\wedge \omega^{\sharp},\, \theta_i\wedge \theta_j\big\rangle,\qquad 1\le i<j\le d,
\]
where $\{\theta_i\}_{i=1}^d$ is the standard orthonormal basis of $\mathbb R^d$.

\begin{theorem}[Ambrozio--Carlotto--Sharp, 2018 \cite{Ambrozio2018DiffGeo}]
Let $(\mathcal N^{n+1},g)$ be a Riemannian manifold that is isometrically embedded in some Euclidean space $\mathbb R^d$. Let $M^n$ be a closed embedded minimal hypersurface of $(\mathcal N^{n+1},g)$. Assume that for every nonzero vector field $X$ on $M^n$,
\begin{equation}\label{eq:ACS_pinch}
\begin{split}
&\int_M \bigl[\mathrm{tr}_M\bigl(\mathrm{Rm}^{\mathcal N}(\cdot,X,\cdot,X)\bigr) + \mathrm{Ric}^{\mathcal N}(N,N)|X|^2\bigr]\,dM \\
&\quad > \int_M \bigl[|\mathit{I\!I}(\cdot,X)|^2 - |\mathit{I\!I}(X,N)|^2 + \bigl(|\mathit{I\!I}(\cdot,N)|^2 - |\mathit{I\!I}(N,N)|^2\bigr)|X|^2\bigr]\,dM,
\end{split}
\end{equation}
where $\mathrm{Rm}^{\mathcal N}$ denotes the Riemann curvature tensor of $\mathcal N^{n+1}$, $\mathit{I\!I}$ denotes the second fundamental form of $\mathcal N^{n+1}$ in $\mathbb R^d$, and $N$ is a local unit normal vector field on $M^n$. Then
\begin{equation}
\mathrm{Index}(M)\;\ge\;\frac{2}{d(d-1)}\, b_1(M).
\end{equation}
\end{theorem}

This gives a \emph{linear-growth} lower bound consistent with the conjectural Schoen--Marques--Neves picture. Moreover, Ambrozio--Carlotto--Sharp verify in \cite{Ambrozio2018DiffGeo} that their pinching condition holds in compact rank one symmetric spaces (CROSS) and other standard positively curved ambient manifolds, yielding a uniform linear lower bound for $\mathrm{Index}$ in terms of $b_1$, with the constant depending only on the Euclidean embedding dimension $d$.

Ambrozio, Carlotto, and Sharp then extended the method to the free boundary setting \cite{Ambrozio2018MathAnn}, studying compact free boundary minimal hypersurfaces in strictly mean-convex domains $\Omega\subset\mathbb R^{n+1}$. Here the role of $b_1(M)$ is played by the relative homology dimension $\dim H_1(M,\partial M;\mathbb R)$. This uses Hodge theory on manifolds with boundary and involves tangential and normal harmonic forms satisfying suitable boundary conditions. In the free boundary setting, integration by parts in the Bochner formula produces additional boundary terms. To overcome this, one exploits the Hodge correspondence between relative homology and harmonic forms that are tangential or normal along the boundary.

\begin{theorem}[Ambrozio, Carlotto \& Sharp, 2018 \cite{Ambrozio2018MathAnn}]
Let $\Omega\subset \mathbb R^{n+1}$ be a strictly mean-convex domain, and let $n\ge 2$. If $M^n$ is a compact, orientable, properly embedded free boundary minimal hypersurface in $\Omega$, then
\begin{equation}
\mathrm{Index}(M)\;\ge\;\frac{2}{n(n+1)}\, \dim H_1(M,\partial M;\mathbb R).
\end{equation}
\end{theorem}

For a free boundary minimal surface $M^2$ in $\mathbb R^3$ of genus $g$ with $r$ boundary components, one has
\[
\dim H_1(M,\partial M;\mathbb R)=2g+r-1,
\]
hence
\[
\mathrm{Index}(M)\;\ge\;\frac{1}{3}\,(2g+r-1)\qquad\text{\cite{Ambrozio2018MathAnn}.}
\]
It is worth noting that in the surface case $n=2$, Ambrozio--Carlotto--Sharp actually prove this estimate under the weaker assumption that $\Omega$ is merely mean convex (not necessarily strictly mean convex), by exploiting a result of Ros \cite{Ambrozio2018MathAnn}.
In particular, \cite{Ambrozio2018MathAnn} points out that the family of genus-zero free boundary minimal surfaces with arbitrarily many boundary components constructed by Fraser and Schoen (2016) in the unit ball has Morse index growing linearly with the number $r$ of boundary components.

In recent years, with the rise of the geometry of weighted manifolds---especially through the study of self-shrinkers for mean curvature flow---the method based on ambient embeddings and harmonic forms has been extended to $f$-minimal hypersurfaces. Impera, Rimoldi, and Savo \cite{Impera2020Rev}, working in the weighted Euclidean space $(\mathbb R^{m+1},g_{\mathrm{can}},e^{-f}dV)$, used the weighted Hodge $f$-Laplacian and the weighted stability operator $L_f$ to prove that compact self-shrinkers (corresponding to $f(x)=|x|^2/2$) satisfy the affine lower bound
\begin{equation}
\mathrm{Ind}_f(\Sigma)\;\ge\;\frac{2}{m(m+1)}\,b_1(\Sigma)\;+\;m+1.
\end{equation}

Impera and Rimoldi \cite{Impera2020Ann} then generalized the Ambrozio--Carlotto--Sharp construction to more general weighted ambient manifolds. In the weighted setting, one replaces ordinary harmonic $1$-forms by $f$-harmonic $1$-forms, characterized by $d\omega=0$ and $\delta_f\omega=0$, and modifies the curvature hypothesis into an extrinsic pinching condition involving the Bakry--\'Emery Ricci tensor
\[
Ric_f=\mathrm{Ric}+\mathrm{Hess}\,f.
\]
In this way one again obtains linear control of $\mathrm{Ind}_f$ by $b_1$ in a broad class of weighted ambient spaces.

Chen, Ge, and Zhang (2023) \cite{Chen2023} further studied compact free boundary $f$-minimal hypersurfaces in weighted manifolds $(\mathcal N^{n+1},g,e^{-f})$, giving a unified framework that simultaneously covers the closed case, the Euclidean free boundary case, and more general weighted ambient spaces. Here $\mathrm{Ric}^{\mathcal N}_f = \mathrm{Ric}^{\mathcal N}+\mathrm{Hess}^{\mathcal N} f$ denotes the Bakry--\'Emery Ricci tensor.

\begin{theorem}[Chen, Ge \& Zhang, 2023 \cite{Chen2023}]\label{thm:CGZ_general}
Let $M^n$ be a compact, orientable, free boundary $f$-minimal hypersurface of a weighted Riemannian manifold $(\mathcal N^{n+1},g,e^{-f}dV_{\mathcal N})$. Assume that $\mathcal N^{n+1}$ is isometrically immersed in some Euclidean space $\mathbb R^d$.
\begin{enumerate}
\item If for every nonzero tangential $f$-harmonic $1$-form $\omega\in\mathcal H^1_{Nf}(M)$,
\begin{equation*}
\begin{split}
&\int_M \bigl[\mathrm{Ric}^{\mathcal N}_f(\omega^\sharp,\omega^\sharp) + \mathrm{Ric}^{\mathcal N}_f(N,N)|\omega|^2 - K^{\mathcal N}(\omega^\sharp,N)\bigr]e^{-f}\,dV_M \\
&\quad > \int_M \sum_{k=1}^n \bigl[|\mathit{I\!I}(e_k,\omega^\sharp)|^2 + |\mathit{I\!I}(e_k,N)|^2|\omega|^2\bigr]e^{-f}\,dV_M \\
&\qquad - \int_{\partial M}\bigl[\mathit{I\!I}^{\partial\mathcal N}(\omega^\sharp,\omega^\sharp) + \mathit{I\!I}^{\partial\mathcal N}(N,N)|\omega|^2\bigr]e^{-f}\,dV_{\partial M},
\end{split}
\end{equation*}
then $\;\mathrm{Index}_f(M)\;\ge\;\frac{2}{d(d-1)}\,\dim H_1(M;\mathbb R)$.
\item If the analogous inequality holds for every nonzero normal $f$-harmonic $1$-form $\omega\in\mathcal H^1_{Tf}(M)$, with $\mathit{I\!I}^{\partial\mathcal N}(\omega^\sharp,\omega^\sharp)+\mathit{I\!I}^{\partial\mathcal N}(N,N)|\omega|^2$ on the boundary replaced by $H^{\partial\mathcal N}_f|\omega|^2$, then $\;\mathrm{Index}_f(M)\;\ge\;\frac{2}{d(d-1)}\,\dim H_{n-1}(M;\mathbb R)$.
\end{enumerate}
Here $K^{\mathcal N}$ is the sectional curvature of $\mathcal N$, $\mathit{I\!I}$ is the second fundamental form of $\mathcal N^{n+1}$ in $\mathbb R^d$, $\mathit{I\!I}^{\partial\mathcal N}$ is the scalar second fundamental form of $\partial\mathcal N$ in $\mathcal N$ with respect to the inward unit normal, $H^{\partial\mathcal N}_f=H^{\partial\mathcal N}+\langle\nabla f,\nu\rangle$ is the $f$-mean curvature of $\partial\mathcal N$, $N$ is a unit normal of $M$ in $\mathcal N$, and $\{e_1,\ldots,e_n\}$ is a local orthonormal frame on $M$.
\end{theorem}

In the Euclidean-domain setting, the above integral conditions can be replaced by simpler pointwise hypotheses.

\begin{theorem}[Chen, Ge \& Zhang, 2023 \cite{Chen2023}]\label{thm:CGZ_Euclidean}
Let $\mathcal N^{n+1}$ be a compact domain in $\mathbb R^{n+1}$, and let $M^n$ be a compact, orientable, free boundary $f$-minimal hypersurface of the weighted manifold $(\mathcal N^{n+1},g_{\mathrm{can}},e^{-f}dV_{\mathcal N})$ with $\mathrm{Hess}_{\mathcal N} f\ge 0$. Then:
\begin{enumerate}
\item if $\partial\mathcal N$ is strictly two-convex in $\mathcal N$, then
\[
\mathrm{Index}_f(M)\;\ge\;\frac{2}{n(n+1)}\,\dim H_1(M;\mathbb R);
\]
\item if $\partial\mathcal N$ is strictly $f$-mean-convex in $\mathcal N$, then
\[
\mathrm{Index}_f(M)\;\ge\;\frac{2}{n(n+1)}\,\dim H_{n-1}(M;\mathbb R).
\]
\end{enumerate}
\end{theorem}

Compared with earlier work, \cite{Chen2023} combines the $f$-Bochner formula with precise control of weighted boundary terms and incorporates, through strict two-convexity or strict $f$-mean-convexity of the boundary, the closed/free-boundary, unweighted/weighted, and Euclidean/general ambient situations into a single extrinsic pinching framework.

\section{Symmetry Reduction and Isoparametric Foliations: Existence Constructions in Euclidean Space and the Unit Ball}\label{sec:existence_equivariant}

\subsection{Complete minimal hypersurfaces via group actions and isoparametric foliations}

In the geometric analysis of minimal submanifolds, the search for higher-dimensional examples typically faces the formidable analytic difficulty of solving complicated nonlinear systems of partial differential equations. A classical idea is to exploit symmetry---either of the ambient space or of the problem itself---to reduce the high-dimensional problem to a lower-dimensional ordinary differential equation or dynamical system. The foundational work in this direction goes back to the pioneering paper of Hsiang and Lawson in 1971 \cite{HsiangLawson1971}. They systematically studied isometric actions of a compact Lie group $G$ on a Riemannian manifold $M$ and established a celebrated reduction theorem, showing that the problem of finding high-dimensional minimal submanifolds can be reduced to the problem of finding geodesics or lower-dimensional minimal submanifolds in the quotient space. This reduction framework was subsequently applied widely to concrete group actions and geometric classification problems.

Alencar (1993) studied the standard action of $SO(m)\times SO(m)$ on $\mathbb R^{2m}$ \cite{Alencar1993}. By analyzing profile curves in the orbit plane, he obtained classification and topological descriptions of complete minimal hypersurfaces, at least in certain dimensions. He proved that, apart from the minimal quadratic cone, all such minimal hypersurfaces are either strictly embedded and topologically homeomorphic to $\mathbb R^m\times S^{m-1}$, or immersed with infinitely many self-intersections and topologically homeomorphic to $\mathbb R\times S^{m-1}\times S^{m-1}$.

\begin{theorem}[Alencar, 1993 \cite{Alencar1993}]
When $m=2$ or $m=3$, every complete $SO(m)\times SO(m)$-invariant minimal hypersurface $M\subset\mathbb R^{2m}$ is either
\begin{enumerate}
\item embedded and topologically homeomorphic to $\mathbb R^m\times S^{m-1}$, or
\item immersed with infinitely many self-intersections and topologically homeomorphic to $\mathbb R\times S^{m-1}\times S^{m-1}$.
\end{enumerate}
Moreover, in both cases the hypersurface intersects the corresponding minimal quadratic cone outside a compact set and can approach that cone arbitrarily closely.
\end{theorem}

Later, Alencar, Barros, Palmas, and collaborators carried out a detailed classification of minimal hypersurfaces invariant under the action of $O(m)\times O(n)$ \cite{Alencar2005}. This generalizes Alencar's 1993 analysis from the case $m=n$ (mainly for $m=2,3$) to the case of possibly different dimensions:
\[
\mathbb{R}^{m+n} = \mathbb{R}^m \times \mathbb{R}^n, \quad G = O(m) \times O(n), \quad m, n \ge 3.
\]
Under this action, a $G$-invariant hypersurface is again described by a profile curve $\gamma(t)=(x(t),y(t))$ in the orbit space
\[
Q=\{(x,y):x\ge 0,\ y\ge 0\}.
\]
The minimality condition reduces to a second-order ODE for $\gamma$, which can be rewritten as an autonomous planar dynamical system and analyzed in terms of asymptotic behavior, crossing properties, and embedded versus self-intersecting solutions.

When $m+n\le 7$, every $O(m)\times O(n)$-invariant minimal hypersurface belongs to one of three classes: the cone $C_{m,n}$ generated by the ray
\[
y = \sqrt{\frac{n-1}{m-1}}\,x,
\]
self-intersecting immersed hypersurfaces that cross this cone countably many times and are asymptotic to it at infinity, and embedded hypersurfaces that also cross the cone infinitely many times and are asymptotic to it, while meeting one of the coordinate-axis boundaries $\mathbb R^m\times\{0\}$ or $\{0\}\times\mathbb R^n$ orthogonally. From the stability viewpoint, every complete hypersurface of this type has infinite Morse index.

When $m+n\ge 8$, the classification changes to four classes. In addition to the basic cone $C_{m,n}$, there are immersed hypersurfaces asymptotic to the cone from both ends while staying on one side of it, embedded hypersurfaces that cross the cone exactly once and are asymptotic to it at both ends, and embedded hypersurfaces that are asymptotic to the cone on one side and meet a coordinate-axis boundary orthogonally. In dimensions $N=m+n\ge 8$, there therefore exist complete embedded stable minimal hypersurfaces whose topology is not that of $\mathbb R^{N-1}$.

The analysis in \cite{Alencar2005} remains based on the Hsiang--Lawson equivariant reduction \cite{HsiangLawson1971}: first the minimality equation is transformed into an ODE for the profile curve, and then phase-plane methods are used to characterize embeddedness, the number of self-intersections, and asymptotic behavior; stability is determined through second variation and spectral analysis, hence through Morse index considerations.

Wang (1994) observed that the reduction mechanism does not fundamentally depend on group orbits themselves, but rather on the existence of a foliation structure producing a good quotient metric and volume density. He extended the preceding group-orbit reduction theory to the more general setting of isoparametric geometry \cite{Wang1994}. Given an isoparametric family of hypersurfaces in $S^{n-1}$, with number of distinct principal curvatures $g\in\{1,2,3,4,6\}$ and multiplicities $(m_1,m_2)$ satisfying the M\"unzner relation
\[
n - 2 = \frac{m_1 + m_2}{2} g,
\]
there exists a Cartan--M\"unzner polynomial $F$ of degree $g$ whose level sets on $S^{n-1}$ are precisely the leaves of the foliation. Extending this spherical foliation radially to $\mathbb R^n$, one calls a hypersurface \textbf{$F$-invariant} if it is a union of such leaves.

Wang explicitly pointed out that this class contains, as a special case, the cohomogeneity-one group-invariant situations of Hsiang and Lawson. He proved that the problem of finding $F$-invariant minimal hypersurfaces in $\mathbb R^n$ is completely equivalent to finding global geodesics in a two-dimensional sector domain $D_g$ endowed with a degenerate Riemannian metric. After a suitable reparametrization, the problem becomes a nonlinear autonomous dynamical system whose topology and asymptotic behavior depend sharply on the algebraic relation between the ambient dimension $n$ and the number $g$ of principal curvatures. The essential point is that Wang reduces the $F$-invariant minimal hypersurface problem to the study of two-dimensional profile curves that are geodesics for a weighted metric, thereby converting a PDE/variational problem into an ODE/dynamical-systems problem. He developed the full framework of isoparametric maps, profile geodesics, and phase-plane dynamics; classified all complete $F$-invariant minimal hypersurfaces (cone/type I/type II); and analyzed stability and area-minimizing properties.

From this unified perspective, the $O(m)\times O(n)$ action corresponds precisely to the isoparametric case $g=2$. Hence the later classification of minimal hypersurfaces under $O(m)\times O(n)$ by Alencar, Barros, and Palmas \cite{Alencar2005} can essentially be viewed as a direct specialization of Wang's theory to the parameter choice
\[
(g,m_1,m_2)=(2,m-1,n-1).
\]
Wang's work thus provides a global framework unifying many earlier special cases involving rotational symmetry and product-group symmetry.

\subsection{Free boundary minimal hypersurfaces via equivariant reduction in the unit ball}

In the closed setting, one emphasizes completeness and asymptotic behavior at infinity, such as asymptotics to a minimal cone $C(M^\ast)$ and the number of crossings with that cone. In the free boundary setting, the object is restricted to the ball $B^{n+1}$, and the natural boundary condition becomes
\[
\partial\Sigma\subset \partial B^{n+1},\qquad \Sigma\perp \partial B^{n+1}.
\]
After equivariant reduction, this free boundary condition is converted into an orthogonality/radial condition at the endpoints of the profile curve in the two-dimensional parameter domain. In other words, the ODE or dynamical system itself does not change; what changes is the choice of orbit segment. In the closed setting one chooses a trajectory satisfying global completeness or prescribed asymptotics, whereas in the free boundary setting one selects a segment that satisfies the orthogonality endpoint condition at the finite radius $r=1$.

Thus the symmetry-reduction and isoparametric phase-plane methods developed in the whole-space setting can be transferred almost verbatim to compact domains with free boundary conditions.

Freidin, Gulian, and McGrath (2017) were the first to implement this transfer systematically within the group-action framework, constructing free boundary minimal hypersurfaces with $O(m)\times O(n)$ symmetry \cite{Freidin2017}. They showed that the same bifurcation threshold already visible in the closed theory reappears in the unit ball: when $m+n\ge 8$, one obtains a family of embedded examples, whereas when $m+n<8$, one can construct infinitely many embedded examples. This threshold is exactly parallel to the condition $n<4g$ versus $n\ge 4g$ in Wang's framework when $g=2$.

Siffert and Wuzyk (2022) \cite{SiffertWuzyk2022} then observed that this transfer principle is much more general. They carried Wang's isoparametric reduction and dynamical-systems analysis over in full to the free boundary problem in the unit ball. Their point of view is that isoparametric hypersurface families generalize group orbits on the sphere, so that $F$-invariance generalizes invariance under a subgroup. In particular, the case $g=2$ corresponds to Clifford tori, i.e.\ the orbits of $O(m_1+1)\times O(m_2+1)$, so the framework conceptually unifies the earlier $O(m)\times O(n)$-symmetric constructions and provides a rigorous path toward extrapolation.

Siffert and Wuzyk proved that for any given isoparametric triple $(g,m_1,m_2)$, which determines the ambient dimension by
\[
n = \frac{m_1+m_2}{2} g + 2,
\]
one can construct $F$-invariant free boundary minimal hypersurfaces in the unit ball $B^n\subset\mathbb R^n$. More precisely:
\begin{enumerate}
    \item when $n<4g$, corresponding to focal-type behavior in the dynamical system, one can construct for each positive integer $k$ a free boundary minimal hypersurface $\Sigma_{g,m_1,m_2}^k$;
    \item when $n\ge 4g$ and excluding several exceptional triples---the main exceptions are $(2,1,5)$, $(4,1,6)$, and the cases obtained by exchanging $m_1$ and $m_2$---one obtains a new free boundary minimal hypersurface $\Omega_{g,m_1,m_2}$.
\end{enumerate}

They state explicitly in their introduction that the first family $\Sigma_{g,m_1,m_2}^k$ had already appeared in Wang's 1994 analysis, while the second family $\Omega_{g,m_1,m_2}$ is genuinely new.

At a structural level, Siffert--Wuzyk rely precisely on the ingredients of Wang's theory: the isoparametric triple $(g,m_1,m_2)$ coming from a spherical isoparametric foliation, the induced ODE/dynamical system satisfied by the reduced profile curve, and the fact that the free boundary condition can be rewritten as a simple angle condition along that curve.

\end{document}